\documentclass[12pt,a4paper,twocolumn]{article}
\usepackage{amsmath,amssymb,amsfonts,amsthm}
\usepackage{mathrsfs} 
\theoremstyle{definition}
\newtheorem{dfn}{Definition}
\theoremstyle{theorem}
\newtheorem{thm}{Theorema}

\newtheorem{cor}{Corollary}
\newtheorem{prp}{Proposition}

\begin{document}
\begin{center}
{\normalsize \bf State estimation for dynamical system described by linear
  equation with uncertainty\footnote{to appear in Ukrainian mathematical journal, 2009}}\\[5ex]
{\bf Serhiy Zhuk}\\
Faculty of Cybernetics\\
Taras Shevchenko Kyiv National University, Ukraine\\
e-mail: {\it beetle@unicyb.kiev.ua}\\[5ex]
\end{center}
\par\textbf{Abstract.} In this paper we investigate a problem of state
estimation for the dynamical system described by the linear operator equation
with unknown parameters in Hilbert space. We present explicit expressions for
linear minimax estimation and error provided that any pair of uncertain
parameters belongs to the quadratic bounding set. As an application of the main
result we present the solution of minimax estimation problem for the linear
descriptor differential equation with constant matrices. 
\normalsize \small{
\section*{Introduction}
}

\normalsize
One of the major problems in applied mathematics is the problem of state
estimation for dynamical systems described by the linear equations with
uncertain parameters. This problem belongs to the so called ``uncertain inverse
problems''. In the strict sense the problem is described as follows: given
some element (for instance measurements of the system output) $y$
 from some functional space one needs to estimate the expression
 $\ell(\theta)$ provided that $\theta$ obey the equation $g(\theta)=0$. This
 problem is non-trivial if there exists non-unique $\theta$ satisfying the
 equation $g(\theta)=0$ and the equality $y=C(\theta)$ holds. In this case the
 estimation problem may be reformulated as follows: given
 $y=C(\theta),\theta\in\Theta,y\in Y$ one needs to find the estimation
 $\widehat{\ell(\theta)}$ of the expression $\ell(\theta)$ provided that
 $g(\theta)=0$ and $C(\cdot),\ell(\cdot)$ are given functions. Note that this
 problem becomes trivial if the equation $y=C(\theta)$ has the unique solution
 $\hat\theta$. Really, in this case we just set
 $\widehat{\ell(\theta)}:=\ell(\hat\theta)$.   

The estimation problem is said to be linear if $\Theta,Y$ are linear spaces
and $C(\cdot),\ell(\cdot)$ are linear mappings. It is the common case when $$
C(\theta)=H\varphi+D\eta,
g(\theta)=L\varphi+Bf, \eqno(*)
$$ where $\theta=(x,f,\eta)\subset X,F,Y$, $H,D,L,B$ are linear mappings. The
linear estimation problem is said to be uncertain if $D\ne0$, $L$ and $B$ are
non-trivial or if $B=0$ then
$\mathrm{N}(L)=\{\varphi:L\varphi=0\}\ne\{0\}$. Note that the choice of 
solution method depends on the ``type of uncertainty'': if $f,\eta$ denotes
realizations of random elements then it's natural to apply probability
methods. This requires an a priori knowledge of distribution characteristics
of the random elements. In the sequel we assume that there is uncertainty in
$(*)$ if distributions of random elements or some deterministic parameters of
the system are partially unknown. An up to date description of
the state of the art in the theory of uncertain estimation
problems with special $\ell,L,H,B,D$ in special spaces is to be found at
\cite{Nakonechnii2004,Chernousko1994,Milanese1991,Kurzhanski1997,Bakan2003,Kuntsevich1992}.  

Classical theory of uncertain estimation problems
\cite{Nakonechnii2004,Milanese1991,Kurzhanski1997} works well when the linear
mapping $L$ in $(*)$ has bounded inverse mapping. The solution of the linear
uncertain estimation problem for linear equations with Noether operator was
introduced in \cite{Podlipenko2005}. Note that the introduced approach is not
suitable when $\mathrm{dim}(L)=\infty$. This is the case for linear
differential-algebraic (DAE) equations \cite{Zhuk2007}. 

The major contribution of this paper is an state estimation approach for
uncertain equations with linear closed operator in Hilbert space. This
approach is valid when $L$ has non-closed range. It generalizes the classical
theory \cite{Nakonechnii2004,Milanese1991,Kurzhanski1997} to the linear
noncausal uncertain differential-algebraic equations. The application of
presented method to the linear equations with Noether operator was described
in \cite{Zhuk2007a}. 

\emph{Notation}. Set $c(G,\cdot)=\sup\{(z,f),f\in G\}$, let $\delta(\mathcal
G,\cdot)$ denotes the indicator function of $\mathcal G$, set
$\mathrm{dom}f=\{x\in\mathcal H:f(x)<\infty\}$,
$f^*(x^*)=\sup_{x}\{(x^*,x)-f(x)\}$, $(L^*c)(u)=\inf\{c(G,z),L^*z=u\}$, $(f
L)(x)=f(Lx)$, $(L^*c)(u)=\inf\{c(G,z),L^*z=u\}$, let $
\mathrm{cl}f=f^{**}$ denotes the closure of $f$ in the Young-Fenhel sense,
$\mathrm{Arginf}_uf(u)$ denotes the set of minimum points of $f$, $P_{L^*}$
denotes the orthogonal projector onto $R(L^*)$, $\partial
f(x)$ denotes the subdifferential of $f$ at $x$ and $(\cdot,\cdot)$ denotes
the inner product in Hilbert space. 
\normalsize \small{
\section*{\bf Problem formulation and definitions.   
}
\normalsize
Suppose that $L\varphi\in\mathscr{G}$ and 
\begin{equation}
  \label{eq:y}
  y=H\varphi+\eta
\end{equation}
The mappings $L,H$ and the set $\mathscr G$ are supposed to be given. The
element $\eta$ is uncertain. Our aim is to solve the inverse problem: to
construct the operator mapping the given $y$ into the estimation
$\widehat{\ell(\varphi)}$ of expression $\ell(\varphi)$ and to calculate the
estimation error $\sigma$. Now let us introduce some definitions. 

The operator $L:\mathcal{H}\mapsto\mathcal{F} $ is assumed to be closed. Its
domain $\mathscr{D}(L)$ is supposed to be a dense subset of the Hilbert space
$\mathcal{H}$, $H\in\mathscr{L}(\mathcal{H},\mathcal{Y})$. Note that the
condition $L\varphi\in\mathscr{G}$ is equal to the following 
\begin{equation}
  \label{eq:Lfi}
  L\varphi=f,
\end{equation}
where $f$ is uncertain and belongs to the given subset $\mathscr{G}$ of the
Hilbert space $\mathcal{F}$. In the sequel $\eta$ is supposed to be a random
$\mathcal{Y}$-valued vector with zero mean so that its correlation $R_\eta\in
\mathscr{R}$,  where $\mathscr{R}$ is some subset
of$\mathscr{L}(\mathcal{Y},\mathcal{Y})$. Also we deal with deterministic
$\eta$ so that  $(f,\eta)\in\mathcal G$, where $\mathcal G$ is some subset of
$\mathcal F\times\mathcal Y$. Note that the realization of $y$ depends $\eta$,
$H$ and $f$. Also it depends 
on elements of $N(L)=\{\varphi\in\mathscr D(L):L\varphi=0\}$ so that
$y=H(\varphi_0+\varphi)+\eta$, where $\varphi_0$ may be thought as inner noise
in the state model \eqref{eq:Lfi}. 

Let $\ell(\varphi)=(\ell,\varphi)$, $\widehat{\ell(\varphi)}=(u,y)+c$. Since
$L,H$ are not supposed to have a bounded inverse mappings the $\ell(\varphi)$
and $\widehat{\ell(\varphi)}$ are not stable with a respect to small
deviations in $f,\eta$. Also $f,\eta$ are supposed to be uncertain. Therefore
we use the minimax design in order to construct the estimation. 
\begin{dfn}\label{ozn1}
  The function $\widehat{\widehat{\ell(\varphi)}}=(\hat{u},\cdot)+\hat{c}$ is
  called the \emph{a priori minimax mean-squared estimation} iff
  $\sigma(\ell,\hat u)=\inf_{u,c}\sigma(\ell,u)$ where
\begin{equation}
    \label{eq:huc}
        \sigma(\ell,u):=\sup_{L\varphi\in \mathscr{G},R_\eta\in\mathscr R}
    M(\ell(\varphi)-\widehat{\widehat{\ell(\varphi)}})^2
  \end{equation}
The number $ \hat{\sigma}(\ell)=\sigma^\frac
12(\ell,\hat u)$ is said to be \emph{the minimax mean-squared error} in the
direction $\ell$. 
\end{dfn}
On the contrast the a posteriori estimation describes the evolution of the central point of the system reachability set $$
(L\varphi,y-H\varphi)\in\mathcal G
$$ consistent with measured output
$y$~\cite{Bertsekas1971,Tempo1985,Nakonechnii1978}. Note that the condition 
$(L\varphi,y-H\varphi)\in\mathcal G$ holds if $\|y\|<C$ for some real $C$. But
 it doesn't hold in our assumptions if $\eta$ is random since $\|R_\eta\|< c$
 doesn't imply $\|y\|<C$ for realizations of $\eta$. Therefore $\eta$ is
 supposed to be deterministic. 
\begin{dfn}\label{ozn2}
The set $$
\mathcal X_y=\{\varphi\in\mathscr D(L):(L\varphi,y-H\varphi)\in\mathcal G\}
$$ is called an a posteriori set. The vector $\hat{\varphi}$ is said to be
minimax a posteriori estimation of $\varphi$ in the direction $\ell$ iff $$
\hat{d}(\ell):=\inf_{\varphi\in\mathcal X_y}\sup_{\psi\in\mathcal X_y}
|(\ell,\varphi)-(\ell,\psi)|=\sup_{\psi\in\mathcal X_y}
|(\ell,\hat{\varphi})-(\ell,\psi)|
$$ The expression $\hat{d}(\ell)$ is called the minimax a posteriori error in
the direction $\ell$.  
\end{dfn}
In the sequel the minimax mean-squared a priori estimation (error) is referred
as minimax estimation (error). 
\normalsize \small{
\section*{
\textbf{Main results}
}

\normalsize
\begin{prp}\label{t:1}
Assume that $\mathscr G$, $\mathscr{R}$ are convex bounded closed subsets of
$\mathcal F$, $\mathscr{L}(\mathcal{Y},\mathcal{Y})$ respectively. For the
given $\ell\in\mathcal H$ the minimax error $\hat\sigma(\ell)$ is finite iff 
\begin{equation}
  \label{eq:setUl}
  \ell-H^*u\in\mathrm{dom}\,\mathrm{cl}(L^*c)\cap
(-1)\mathrm{dom}\,\mathrm{cl}(L^*c)
\end{equation} for some $u\in\mathcal Y$. Under this condition 
  \begin{equation} 
    \label{eq:err:amap:thr}
    \begin{split}
          &\sigma(\ell,u)=\sup_{R_\eta\in\mathscr R}(R_\eta u,u)+\\
          &\frac 14[\mathrm{cl}(L^*c)(\ell-H^*u)+
    \mathrm{cl}(L^*c)(-\ell+H^*u)]^2
 \end{split}
    \end{equation}
where $$
R(L^*)\subset\mathrm{dom}\,\mathrm{cl}(L^*c)\subset\overline{R(L^*)}
$$
If $\mathrm{Arginf}_u\sigma(\ell,u)\ne\varnothing$, then
$\widehat{\widehat{\ell(\varphi)}}=(\hat{u},y)+\hat{c}$, where $$
\hat{u}\in\mathrm{Arginf}_u\sigma(\ell,u)$$
and $$
\hat{c}=\frac 12(\mathrm{cl}(L^*c)(\ell-H^*\hat u)-
\mathrm{cl}(L^*c)(-\ell+H^*\hat u))
$$
\end{prp}
\begin{thm}\label{t:RLintG}
  Suppose that $\mathscr G_1$ is convex bounded closed balanced set and its
   $0\in\mathrm{int}\,\mathscr G_1$. Also assume that 
  $$
  \eta\in\{\eta:M(\eta,\eta)\le 1\}
  $$ Then for the given $\ell\in\mathcal H$ the minimax estimation
  $\hat\sigma(\ell)$ is finite iff $\ell-H^*u\in R(L^*)$  
  for some $u\in\mathcal Y$. Under this condition there exists a unique
  minimax estimation $\hat u\in\mathcal U_\ell$ and 
\begin{equation}
  \label{eq:aprer:amap:thr}
  \begin{split}
    &\sigma(\ell,\hat u)=\min_u\sigma(\ell,u),\\
    &\sigma(\ell,u)=(u,u)+\min_z\{c^2(\mathscr G_1,z),L^*z=\ell-H^*u\}
  \end{split}
\end{equation}
  If $R(L),H(N(L))$ are closed sets then $\hat u$ is determined by the
  following conditions
  \begin{equation}
    \label{eq:Hp0}
    \begin{split}
    &\hat u-Hp_0\in H(\partial I_2(H^*\hat u)),Lp_0=0,\\
    &I_2(w)=\min_z\{c^2(\mathscr G_1,z),L^*z=P_{L^*}(\ell-w)\},
    \end{split}
  \end{equation}
\end{thm}
\begin{cor}\label{t:2}
Let $$
\mathscr G_1=\{f\in\mathcal F:(f,f)\le 1\},
\eta\in\{\eta:M(\eta,\eta)\le 1\},
$$ and suppose that 
\begin{itemize}
\item [1)] $R(L),H(N(L))$ are closed sets;
\item [2)] $R(T)=\{[Lx,Hx],x\in\mathscr D(L)\}$ is closed set.
\end{itemize}
Then only for $\ell\in R(L^*)+R(H^*)$ the unique minimax estimation $\hat u$
is given by $\hat u=H\hat p$, where $\hat p$ is any solution of the equations 
\begin{equation}
  \label{eq:qeuler}
  \begin{split}
    &L^*\hat z=\ell-H^*H\hat p,\\
    &L\hat p=\hat z
  \end{split}
\end{equation}
The minimax error is given by the following expression
\begin{equation*}
  \hat\sigma(\ell)=(\ell,\hat p)^\frac 12
\end{equation*}
\end{cor}
\begin{cor}\label{n:rozveuler}
Assume that linear mappings $
  L:\mathcal H\mapsto\mathcal F$, $
  H\in\mathscr L(\mathcal H,\mathcal Y)$ obey 1) or 2) (Cor.~\ref{t:2}).  
Then~\eqref{eq:qeuler} has a solution $
\hat z\in\mathscr D(L^*),\hat p\in\mathscr D(L)$ 
 iff $\ell=L^*z+H^*u$ for some $
z\in\mathscr D(L^*),u\in\mathcal Y$.
\end{cor}
\begin{cor}\label{n:alt}
Under the conditions of Cor.~\ref{t:2} for any $\ell\in R(L^*)+R(H^*)$ and
some realization of $y(\cdot)$ we have $(\hat u,y)=
(\ell,\hat\varphi)$, where $\hat\varphi$ obeys the equation
\begin{equation}
  \label{eq:alt}
  \begin{split}
    &L^*\hat q=H^*(y-H\hat\varphi),\\
    &L\hat\varphi=\hat q
  \end{split}
\end{equation}
\end{cor}
Consider an a posteriori estimation. 
\begin{prp}\label{t:apo:1}
Let $\mathcal G$ be a convex closed bounded subset of $\mathcal
  Y\times\mathcal F$. The minimax a posteriori error in the direction $\ell$
  is finite iff $\ell\in
  \mathrm{dom}\,c(\mathcal X_y,\cdot)\cap(-1) 
  \mathrm{dom}\,c(\mathcal X_y,\cdot)$ and
  \begin{equation}
    \label{eq:setL}
    \begin{split}
      &R(L^*)+R(H^*)\subset\mathrm{dom}\,c(\mathcal X_y,\cdot)\cap
    (-1)\mathrm{dom}\,c(\mathcal X_y,\cdot)\subset\\
    &\overline{R(L^*)+R(H^*)}
    \end{split}
  \end{equation}
  Under this condition
  \begin{equation}
    \label{eq:apo:esterr}
    \begin{split}
      &(\ell,\hat\varphi)=\frac 12(c(\mathcal X_y,\ell)-c(\mathcal
    X_y,-\ell)),\\
    &\hat{d}(\ell)=\frac 12(c(\mathcal X_y,\ell)+c(\mathcal X_y,-\ell))
    \end{split}
  \end{equation}
\end{prp}
\begin{thm}\label{t:apo:2}
Let $$
\mathcal G=\{(f,\eta):\|f\|^2+\|\eta\|^2\le 1\},
$$ and assume that one of the conditions of Cor.~\ref{t:2} is fulfilled. The
minimax a posteriori estimation $\hat{\varphi}$ obeys the equation 
\begin{equation}
  \label{eq:apoest:amap:thr}
  \begin{split}
    &L^*\hat q=H^*(y-H\hat\varphi),\\
    &L\hat\varphi=\hat q
  \end{split}
\end{equation}
iff $\ell\in R(L^*)+R(H^*)$. The estimation error is given by 
\begin{equation}
  \label{eq:apos-err}
  \hat{d}(\ell)=(1-(y,y-H\hat\varphi))^\frac 12 \hat\sigma(\ell)
\end{equation}
\end{thm}
\begin{cor}\label{amap:thr:vest}
  Assume that the conditions of Theorem~\ref{t:apo:2} are fulfilled and
  $\widehat{\ell(\varphi)}=(\ell,\hat\varphi)$ for any $\ell$, where
  $\hat{\varphi}$ obeys~\eqref{eq:apoest:amap:thr}. Then $\hat\varphi$ gives
  the minimax a posteriori estimation of $\varphi$ so that 
\begin{equation*}
  \begin{split}
    &\inf_{\varphi\in\mathcal X_y}\sup_{x\in\mathcal X_y}\|\varphi-x\|=\\
    &\sup_{x\in\mathcal X_y}\|\hat\varphi-x\|=(1-(y,y-H\hat\varphi)^\frac 12
    \max_{\|\ell\|=1}\hat{\sigma}(\ell)
  \end{split}
\end{equation*}
\end{cor}

Now we shall apply Cor.~\ref{amap:thr:vest} to the linear uncertain estimation
problem for differential-algebraic equation. Taking into account that any DAE
with constant matrices has the SVD canonical form \cite{Bender1987} we assume
without loss of generality that $$
F=\bigl(
\begin{smallmatrix}
  E&&0\\0&&0
\end{smallmatrix}
\bigr),
C=\bigl(
\begin{smallmatrix}
  C_1&&C_2\\C_3&&C_4
\end{smallmatrix}
\bigr)
$$ 
\begin{prp}\label{ade:est:tv1apo}
   Let $t\mapsto x(t)\in\mathbb R^n$ be the solution of $$
\dfrac d{dt}Fx(t)-Cx(t)=f(t),Fx(t_0)=0,
$$ and set $$
\mathcal G=\{(f,\eta):\int_{t_0}^T(\|f(t)\|^2+\|\eta(t)\|^2)dt\le1\}
$$ Then the minimax a posteriori estimation of the function $x(\cdot)$ is
given by $\hat
x(\cdot)$, where $\hat x(t)=[x_1(t),x_2(t)]$, 
\begin{equation}\label{eq:x12q12:mnmx}
  \begin{split}
    &\dot x_1(t)=(C_1-C_2(E+C_4'C_4)^{-1}C'_4C_3)x_1(t)+\\
    &(C_2(E+C_4'C_4)^{-1}C'_2+E)q_1(t)+\\
    &C_2(E+C_4'C_4)^{-1}y_2(t), \\
    &\dot q_1(t)=(-C'_1+C'_3C_4(E+C_4'C_4)^{-1}C'_2)q_1(t)+\\
    &C'_3C_4(E+C_4'C_4)^{-1}y_2(t)-y_1(t)+\\
    &(C'_3(E-C_4(E+C_4'C_4)^{-1}C'_4)C_3+E)x_1(t),\\
    &q_1(T)=0,x_1(t_0)=0,\\
    &x_2(t)=-(E+C_4'C_4)^{-1}C'_4C_3x_1(t)+\\
    &(E+C_4'C_4)^{-1}(C'_2q_1(t)+y_2(t)),\\
    &q_2(t)=-(E-C_4(E+C_4'C_4)^{-1}C'_4)C_3x_1(t)-\\
    &C_4(E+C_4'C_4)^{-1}(C'_2q_1(t)+y_2(t))
  \end{split}
\end{equation}
provided that $y(t)=x(t)+\eta(t),t_0\le t\le T$. The minimax error is given
by $$ 
\sup_{\mathcal X_y}\|x-\hat x\|=(1-\int_{t_0}^T(y,y-\hat x)dt)^\frac 12
\max_{\|\ell\|=1}(\int_{t_0}^T(\ell,p)dt)^\frac 12
$$ where $p(\cdot)$ obeys~\eqref{eq:x12q12:mnmx} with $y(t)=\ell(t)$.  
\end{prp}

\end{document}